\documentclass[11pt]{amsart}

\usepackage{xcolor}
\usepackage{amsmath,amssymb,graphicx,microtype,tikz}
\usetikzlibrary{arrows.meta,calc,decorations.markings}
\usepackage[letterpaper,hmargin=1.20in,top=.90in,bottom=1.00in,footskip=20pt]{geometry}
\usepackage[hidelinks,pdfusetitle]{hyperref}
\hypersetup{pdfauthor={R. İnanç Baykur}}

\newtheorem{mainthm}{Theorem}

\newtheorem{prop}{Proposition}[section]

\newtheorem{lemma}[prop]{Lemma}
\newtheorem{corollary}{Corollary}
\theoremstyle{remark}
\newtheorem{remark}[prop]{Remark}

\newcommand{\x}{\times}
\newcommand{\Z}{\mathbb{Z}}
\newcommand{\R}{\mathbb{R}}
\newcommand{\CP}{{\mathbb{CP}}{}^{2}}
\newcommand{\CPb}{\overline{\mathbb{CP}}{}^{2}}
\newcommand{\eu}{\mathrm{e}}
\newcommand{\Int}{\operatorname{int}}
\newcommand{\gen}[1]{\langle #1\rangle}

\allowdisplaybreaks[2]
\title[Small exotic 4-manifolds with signature zero]{Small exotic 4-manifolds with signature zero}
\author[R. \.{I}. Baykur]{R. \.{I}nan\c{c} Baykur}
\address{Department of Mathematics and Statistics, University of Massachusetts,
  Amherst, MA 01003, USA}
\email{inanc.baykur@umass.edu}
\date{}

\begin{document}

\begin{abstract}
We construct infinitely many pairwise nondiffeomorphic irreducible symplectic \linebreak $4$--manifolds in the homeomorphism classes of
$\#_3(S^2\x S^2)$ and $\#_3(\CP\#\CPb)$. We extend these results to a larger number of summands and other finite cyclic fundamental groups. We also construct complex examples in the homeomorphism classes of $\#_{23}(S^2\x S^2)$ and $\#_{23}(\CP\#\CPb)$.
\end{abstract}

\AddToHookNext{cmd/@settitle/before}{\topskip=12pt\relax}
\maketitle
\begingroup
\linespread{1.06}\selectfont
\vspace{.5\baselineskip}

\section{Introduction}

Signature zero has been a particularly difficult case to conquer in the
geography of small simply connected exotic four-manifolds~\cite[Problem~4.3(b)]{K3}.
Since the first examples of infinitely many exotic smooth structures on $\#_m(S^2\x S^2)$
and $\#_n(\CP\#\CPb)$ for very large odd $m$ and $n$ by J. Park almost 30 years ago \cite{Park98,Park},
the number of summands has been brought down gradually in successive works of Akhmedov, D.~Park,
Sakall{\i} and Yeung~\cite{APspin,APS,ASY}, and most recently to $m=11$ and $n=9$ by Hamada and the author~\cite{BHfib,BH}. The goal of this article is to show that three summands suffice.

\begin{mainthm}\label{thm:main}
Both $\#_3(S^2\x S^2)$ and $\#_3(\CP\#\CPb)$ admit infinitely
many pairwise nondiffeomorphic irreducible smooth structures, infinitely many of them supporting symplectic structures, and infinitely many supporting none.
\end{mainthm}

\noindent In Corollary~\ref{cor:larger}, we extend this to all odd $m,n>1$, and in Corollary~\ref{cor:cyclic} and Remark~\ref{rem:cyclic-extensions} we treat the signature zero geography with finite cyclic fundamental groups. An equivariant version of our constructions yields irreducible exotic examples with $\pi_1=\Z_2$ and $b_2^+=b_2^-=1$ as free quotients (Remark~\ref{rem:free}).

We obtain our small examples much as the small exotica with negative
signature were obtained in~\cite{ABP,AP,APodd,BK,FPS}, using
symplectic fiber sums of products of surfaces and Luttinger surgeries along Lagrangian tori; see~\cite{Thesis,ABP}
for the underlying construction scheme. For some of the simplest choices of models, the
attempt to kill $\pi_1$ is doomed, since these models are and do remain aspherical after torus surgeries~\cite{Fake}, and are therefore not homeomorphic to a simply connected $4$--manifold, such as $\#_m(S^2\x S^2)$
or $\#_n(\CP\#\CPb)$. This led us to analyze perhaps what we may claim to be the next simplest model for signature zero exotica that passes the asphericity test, the symplectic fiber sum of
$\Sigma_2\x\Sigma_2$ along its diagonal with $T^4$ along its
smoothed coordinate tori. We observe that this model is nonaspherical despite both closed summands being aspherical,
and that a nontrivial commutator in the fundamental group of the gluing boundary
dies on the $T^4$ side
(Remark~\ref{rem:model}) --- unlike in the aspherical fake $S^2 \x S^2$s \cite{Fake}. In Section~\ref{sec:construction} we describe
the torus surgeries that kill $H_1$ (Lemma~\ref{lem:homology}) in this model, and in
Section~\ref{sec:groups} we use this key boundary relation to show that the
surgeries also kill $\pi_1$. In Section~\ref{sec:families} we use knot
surgery along a symplectic torus with a simply connected complement we identify in each case, to obtain the infinite families.

\smallskip
Simply connected signature zero compact \emph{complex} surfaces of general type were constructed even earlier by Moishezon, Teicher, Robb, and Z.~Chen~\cite{MT86,MT87,MRT96,ZChen87}, with fairly large topology.
For instance, the four spin examples of signature zero in~\cite{MRT96}
all have  $b_2^+ \approx 8.71 \x 10^{77}$\cite[Section~4]{LeBrunLi}. 
Our second result gives the smallest such examples we are aware of in either parity, with $b_2^+=23$.

\begin{mainthm}\label{thm:complex}
Both $\#_{23}(S^2\x S^2)$ and $\#_{23}(\CP\#\CPb)$ admit an exotic copy that is a minimal compact complex surface of general type.
\end{mainthm}

Each example is obtained by resolving six singularities in an orbifold cover of a quotient of a complex fake projective plane (FPP). The two FPPs belong to different classes in Cartwright and Steger's classification~\cite{CSplanes}. The same result was obtained independently by Zolt\'an Szab\'o in a similar fashion (more about this is discussed below); the proof given in Section~\ref{sec:complex} is the author's. Notably, we observe that one cannot derive a complex example with $b_2^+=b_2^-=3$ and $b_1=0$ by this construction from \emph{any} FPP (Remark~\ref{rem:fpp-obstruction}). The contrast leads to the following questions regarding the difference between
symplectic and complex geography. \textit{How small can $b_2^+=b_2^-$ be for a simply connected compact complex
surface of general type in each parity?} \textit{Can either or both homeomorphism types in
Theorem~\ref{thm:main} be realized by a complex surface?} 

\smallskip
\subsection*{Conventions}
We use product orientations and integral homology unless otherwise indicated.
Write $\#_mY$ for the oriented connected sum of $m$ copies of $Y$, $\CPb$
for the complex projective plane $\CP$ with reversed orientation, and $\Sigma_g^b$ for the compact
oriented genus--$g$ surface with $b$ boundary components, where $\Sigma_g=\Sigma_g^0$.
Read paths from left to right, with $g^h=h^{-1}gh$ and $[g,h]=ghg^{-1}h^{-1}$ in the fundamental group.
A $1/n$ Luttinger surgery sends the new meridian to $\mu(\lambda')^n$,
where $\mu$ is the old meridian and $\lambda'$ the Lagrangian push-off
of the surgery curve. For the symplectic fiber sums and Luttinger surgeries we refer to~\cite{Gompf} and~\cite{ADK}.
For the complex constructions, the reader can turn to~\cite{Keum}
for quotients of FPPs and to~\cite[Chapter~III]{BHPV} for the resolution of surface singularities.

\smallskip
\subsection*{AI use and more on the methodology}

In the course of proving our main result, Theorem~\ref{thm:main}, we used
ChatGPT to triangulate the torus complements and probe which torus
surgeries might kill $\pi_1$ through these calculations. This may be a small data set for a computer, but is quite sizable for a hand check; the triangulation of the $\Sigma_2\x\Sigma_2$ piece
alone has $146{,}016$ 4-dimensional simplices. The auxiliary files give the resulting presentations and
calculations for verifying
Proposition~\ref{p:groups} and the exterior groups
in Lemma~\ref{lem:exteriors} by computer.
Most of the author's effort went into providing a more tractable hand
proof of simple connectivity, where the
vanishing of the boundary commutator in Remark~\ref{rem:model} allowed
us to kill the remaining generators.

For our secondary result, Theorem~\ref{thm:complex}, the necessary covering calculations were carried out by AI. Starting with the group presentations
in Cartwright and Steger's catalogue~\cite{CSplanes},
we looked for covers of the FPP quotients whose resolutions would be simply connected and have signature zero. In each fixed degree this is a finite game with many initial possibilities, and getting the desired topology requires a fine balance. The proofs given here were written by the author.

Besides these calculations, ChatGPT assisted with the figures, the bibliography and proofreading. As the customary acknowledgement goes
these days, the author is responsible for all the mathematical arguments and
exposition.

\smallskip
\subsection*{Acknowledgements}
This work was supported in part by NSF grant DMS-2506431 and a Simons
Foundation Travel Grant. The author thanks Zolt\'an Szab\'o for helpful
discussions on the complex signature zero examples. Szab\'o independently
obtained Theorem~\ref{thm:complex} and kindly agreed to inclusion of this result here.

\endgroup

\newpage
\section{The construction}\label{sec:construction}

Let $P, Q \cong \Sigma_2$ be decomposed as $P=P_1\cup_C P_2$ and
$Q=Q_1\cup_C Q_2$, where $P_i, Q_i \cong \Sigma_1^1$. For an explicit model for each $P_i$ and $Q_i$, take a copy of
the square $[-10,10]^2 \subset \R^2$, identify opposite sides by translation, and remove the
open disk $x^2+y^2<2$. In each surface, glue
the boundary circles of its two punctured tori by $(x,y)\mapsto(y,x)$.
These coordinates identify $P_i$ with $Q_i$ and the two separating
circles with each other, so we write $C$ for either copy. Let $a_i,b_i$ be the horizontal and vertical circles on $P_i$, oriented
in the directions of increasing $x$ (east) and $y$ (north), and let $c_i,d_i$ be the
corresponding circles on $Q_i$. Each coordinate loop goes once around its
circle in the indicated direction.

The normal bundle of the diagonal $\Delta\subset P\x Q$ is isomorphic to
its tangent bundle, so $\Delta^2=\eu(\Delta)=-2$, where $\eu$ denotes the Euler characteristic.
Write $T^2=\R^2/(5\Z)^2$, and take the point $\xi=(1,1)$.
In $T^4=T^2\x T^2$, smoothing the coordinate tori $T^2\x\{\xi\}$ and
$\{\xi\}\x T^2$ at $(\xi,\xi)$ gives a genus--$2$ surface $\Sigma'$
of square $2$.
Match the areas of $\Delta$ and $\Sigma'$ and form the symplectic fiber
sum~\cite{Gompf}
\[
               M=(P\x Q)\mathbin{\#}_{\Delta=\Sigma'}T^4.
\]
The first projection identifies $\Delta$ with $P$. Its two punctured
tori are the diagonals $\Delta_i=\{(p,p)\mid p\in P_i\}\subset P_i\x Q_i$, using
$P_i\cong Q_i$. As usual, the smoothing removes a disk from each coordinate torus
and joins the two boundary circles by an annulus. Identify $\Delta_1$
and $\Delta_2$ with the punctured $T^2\x\{\xi\}$ and $\{\xi\}\x T^2$,
respectively, matching the oriented horizontal and vertical circles. Equation~\eqref{eq:flat}
specifies the gluing on the normal meridian and on based lifts of these
circles to the normal circle bundles.

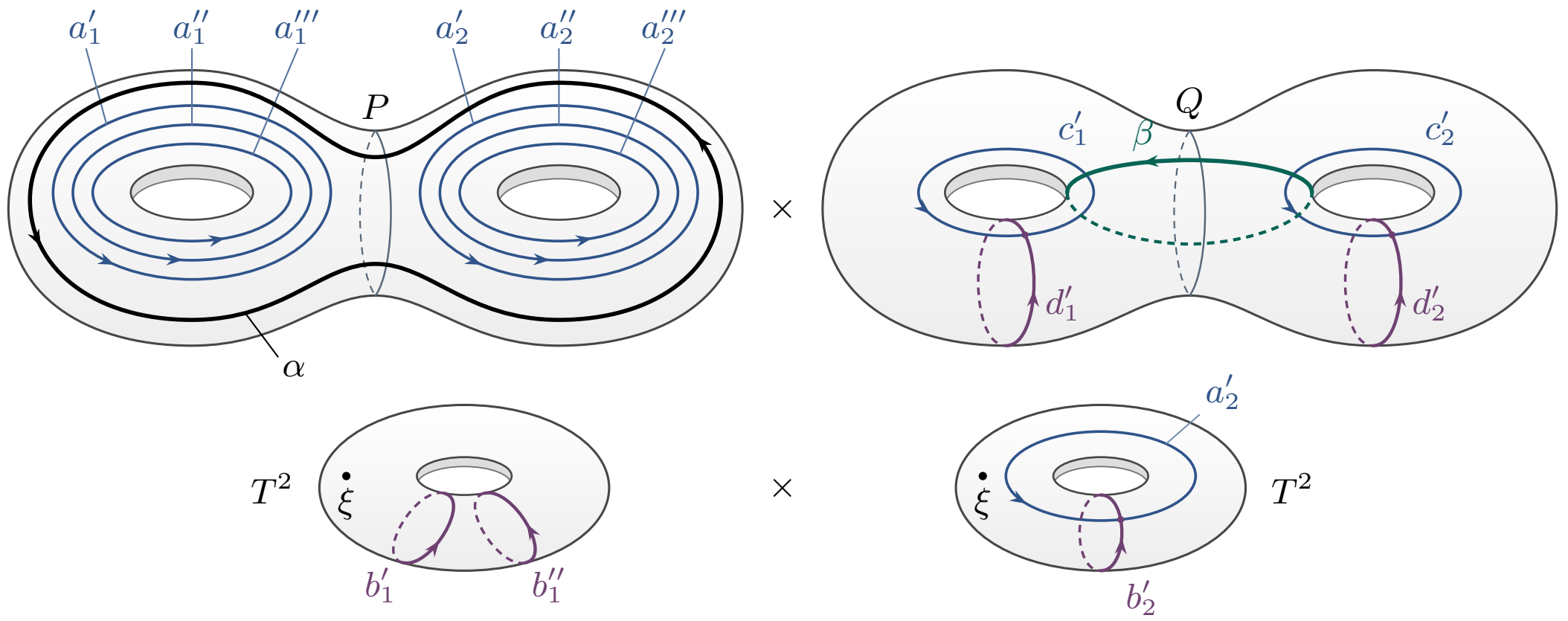
\begin{figure}[ht]
\centering
\begin{tikzpicture}[x=1cm,y=1cm,font=\footnotesize,inner sep=1.3pt,
                    line cap=round,line join=round]
  \definecolor{horizink}{RGB}{48,84,138}
  \definecolor{vertink}{RGB}{111,65,114}
  \colorlet{alphaInk}{black}
  \definecolor{betaInk}{RGB}{7,100,85}
  \definecolor{seamink}{RGB}{87,104,121}
  \tikzset{
    outline/.style={draw=black!72,line width=.58pt},
    surface/.style={outline,top color=white,bottom color=black!7},
    directed/.style={postaction={decorate},decoration={markings,
      mark=at position #1 with {\arrow{Stealth[length=1.45mm,width=1.05mm]}}}},
    parallel/.style={draw=horizink,line width=.68pt},
    meridian/.style={draw=vertink,line width=.83pt},
    back arc/.style={dash pattern=on 1.9pt off 1.45pt,line cap=butt}}
  \def\opening{
    \path[fill=black!13] (0,.13) ellipse[x radius=.37,y radius=.185];
    \begin{scope}
      \clip (0,.13) ellipse[x radius=.37,y radius=.185];
      \path[fill=white,draw=black!55,line width=.35pt]
        (0,.045) ellipse[x radius=.368,y radius=.178];
    \end{scope}
    \draw[outline,line width=.50pt] (0,.13) ellipse[x radius=.37,y radius=.185];}
  \def\genustwo{
    \path[surface] (0,.55)
      .. controls (.32,.55) and (.57,.96) .. (1.11,.96)
      .. controls (1.77,.96) and (2.22,.58) .. (2.22,.02)
      .. controls (2.22,-.54) and (1.78,-.91) .. (1.11,-.91)
      .. controls (.56,-.91) and (.33,-.57) .. (0,-.57)
      .. controls (-.33,-.57) and (-.56,-.91) .. (-1.11,-.91)
      .. controls (-1.78,-.91) and (-2.22,-.54) .. (-2.22,.02)
      .. controls (-2.22,.58) and (-1.77,.96) .. (-1.11,.96)
      .. controls (-.57,.96) and (-.32,.55) .. (0,.55) -- cycle;
    \draw[seamink,back arc,line width=.42pt]
      (0,.55) .. controls (-.125,.36) and (-.125,-.38) .. (0,-.57);
    \draw[seamink,line width=.47pt]
      (0,-.57) .. controls (.125,-.38) and (.125,.36) .. (0,.55);
    \foreach \xc in {-1.11,1.11}{\begin{scope}[shift={(\xc,0)}]\opening\end{scope}}}
  \def\torus{
    \path[surface] (0,.015) ellipse[x radius=1.13,y radius=.81];\opening}
  \def\merid#1#2#3{
    \coordinate (outer) at #1; \coordinate (inner) at #2;
    \draw[vertink,back arc,line width=.65pt] (inner) .. controls ($(inner)-#3$) and ($(outer)-#3$) .. (outer);
    \draw[meridian,directed=.45] (outer)
      .. controls ($(outer)+#3$) and ($(inner)+#3$) .. (inner);}
  \def\longitude#1#2#3{
    \draw[parallel,directed=#3] (#1,.13)
      arc[start angle=0,end angle=360,x radius=#1,y radius=#2] -- cycle;}

  \foreach \cx/\factorname in {-3.62/P,3.62/Q}{
    \begin{scope}[shift={(\cx,0)},x=1.47cm,y=1.31cm]
      \genustwo
      \node[above=1.4pt] at (0,.55) {$\factorname$};
      \foreach \xc/\idx in {-1.11/1,1.11/2}{
        \begin{scope}[shift={(\xc,0)}]
          \ifdim\cx pt<0pt
            \longitude{.84}{.59}{.64}
            \longitude{.72}{.46}{.735}
            \longitude{.60}{.33}{.82}
            \foreach \rx/\ry/\theta/\lx/\primes in
              {.84/.59/128/-.64/1,.72/.46/90/0/2,.60/.33/52/.64/3}{
              \node[text=horizink,anchor=south,inner sep=1pt] (name)
                at (\lx,1.105) {$a_{\idx}^{\ifcase\primes\or\prime\or
                   \prime\prime\or\prime\prime\prime\fi}$};
              \draw[horizink!80,line width=.4pt]
                (name.south) -- ({\rx*cos(\theta)},{.13+\ry*sin(\theta)});
            }
          \else
            \longitude{.53}{.295}{.56}
            \merid{(0,-.91)}{(0,-.055)}{(.225,0)}
            \fill[vertink] (.11457,-.15803) circle[radius=.88pt];
            \node[text=horizink,above] at (.41,.41) {$c_{\idx}'$};
            \node[text=vertink,right] at (.20,-.62) {$d_{\idx}'$};
          \fi
        \end{scope}}
      \ifdim\cx pt<0pt
        \draw[alphaInk,line width=1.05pt,directed=.28,directed=.79]
          (0,.37)
          .. controls (-.33,.37) and (-.53,.875) .. (-1.11,.875)
          .. controls (-1.69,.875) and (-2.09,.52) .. (-2.09,.08)
          .. controls (-2.09,-.39) and (-1.70,-.735) .. (-1.11,-.735)
          .. controls (-.55,-.735) and (-.33,-.355) .. (0,-.355)
          .. controls (.33,-.355) and (.55,-.735) .. (1.11,-.735)
          .. controls (1.70,-.735) and (2.09,-.39) .. (2.09,.08)
          .. controls (2.09,.52) and (1.69,.875) .. (1.11,.875)
          .. controls (.53,.875) and (.33,.37) .. (0,.37) -- cycle;
        \draw[alphaInk,line width=.4pt] (-.59,-.985) -- (-.79,-.685);
        \node[text=alphaInk,anchor=north west,inner sep=.9pt] at (-.59,-.985) {$\alpha$};
      \else
        \draw[betaInk,back arc,line width=.8pt]
          (-.74,.13) .. controls (-.74,-.08) and (-.32,-.22) .. (0,-.22)
          .. controls (.32,-.22) and (.74,-.08) .. (.74,.13);
        \draw[betaInk,line width=1.06pt,directed=.68]
          (.74,.13) .. controls (.74,.28) and (.32,.35) .. (0,.35)
          .. controls (-.32,.35) and (-.74,.28) .. (-.74,.13);
        \node[text=betaInk,above=0pt] at (-.28,.37) {$\beta$};
      \fi
    \end{scope}}
  \node at (0,.02) {$\x$};

  \foreach \cx/\idx in {-2.83/1,2.83/2}{
    \begin{scope}[shift={(\cx,-2.47)},x=1.14cm,y=.91cm]
      \torus
      \ifnum\idx=1
        \merid{(-.48,-.716)}{(-.157,-.038)}{(.22,-.065)}
        \merid{(.48,-.716)}{(.157,-.038)}{(.22,.065)}
        \node[text=vertink,below left=0pt] at (-.48,-.71) {$b_1'$};
        \node[text=vertink,below right=0pt] at (.48,-.71) {$b_1''$};
        \node[left=5pt] at (-1.13,.01) {$T^2$};
      \else
        \longitude{.74}{.435}{.58}
        \merid{(0,-.795)}{(0,-.055)}{(.22,0)}
        \fill[vertink] (.15563,-.29527) circle[radius=.88pt];
        \draw[horizink!72,line width=.35pt] (.77,.73) -- (.51,.445);
        \node[text=horizink,above right=0pt] at (.77,.73) {$a_2'$};
        \node[text=vertink,below right=0pt] at (.15,-.83) {$b_2'$};
        \node[right=5pt] at (1.13,.01) {$T^2$};
      \fi
      \fill (-.92,.13) circle[radius=1.05pt] node[below=1.5pt] {$\xi$};
    \end{scope}}
  \node at (0,-2.46) {$\x$};
\end{tikzpicture}
\caption{On $P$, the black curve $\alpha$ and the outer blue parallels
$a_1',a_2'$ bound a pair of pants. On $Q$, the green curve $\beta$
misses $d_1',d_2'$ and meets each $c_i'$ once, with
$c_1'\cdot\beta=1$ and $c_2'\cdot\beta=-1$. The ninth surgery torus
is $S=\alpha\x\beta$.}

\label{fig:products}
\end{figure}

In the square models for $P_i$ and $Q_i$, the parallels
$a_i',a_i'',a_i'''$ lie at $y=3,4,5$, respectively, while $c_i',d_i'$
lie at $y=-4$, $x=-4$. In the $i$-th factor of $T^4$ we again write
$a_i,b_i$ for the oriented horizontal and vertical circles through $(0,0)$.
Their parallels are shown in Figure~\ref{fig:products}. The circles
$b_1',b_1''$ lie at $x=2,3$, respectively, and $a_2',b_2'$ at
$y=2$, $x=2$.

The following eight Luttinger surgeries give $X$. For each circle choose the normal for which the positive tangent
followed by the normal orients the surface. Orient the normal disk
to the torus by the two normals in factor order, and give the meridian
its boundary orientation. The last column lists
a product dual torus for each surgery torus, up to orientation.
In Lemma~\ref{lem:homology} we choose disjoint parallels for these duals
and show that the meridians of all nine tori vanish in $H_1$ of their
common exterior. Section~\ref{sec:groups} gives the based meridian relations in $\pi_1$.
\begin{equation}\label{eq:surgeries}
\renewcommand{\arraystretch}{1.14}
\setlength{\arraycolsep}{4pt}
\begin{array}{@{}llccc@{}}
  & \underline{\text{\strut surgery torus}} & \underline{\text{\strut curve}}
  & \underline{\text{\strut coefficient}} & \underline{\text{\strut dual torus}} \\[4pt]
  P\x Q & L_{12}=a_1'\x c_2'  & c_2'   & +1 & b_1\x d_2\\
        & L_{21}=a_2'\x c_1'  & c_1'   & +1 & b_2\x d_1\\
        & L'_{12}=a_1''\x d_2' & d_2'   & -1 & b_1\x c_2\\
        & L'_{21}=a_2''\x d_1' & d_1'   & -1 & b_2\x c_1\\
        & T_1=a_1'''\x c_1'   & a_1''' & +1 & b_1\x d_1\\
        & T_2=a_2'''\x c_2'   & a_2''' & +1 & b_2\x d_2\\[2pt]
  T^4   & T'_1=b_1'\x a_2'    & b_1'   & -1 & a_1\x b_2\\
        & T'_2=b_1''\x b_2'   & b_2'   & +1 & a_1\x a_2
\end{array}
\end{equation}

\smallskip
The first circle factors of the six tori in $P\x Q$ are disjoint. Those in $P_i\x Q_j$
with $i\ne j$ miss $\Delta$, and the circles of $T_i$ have distinct heights
$5,-4$ in the identified handle. Sliding the first circle across the annulus
$3\le y\le5$ with $c_i'$ fixed gives an isotopy from $a_i'\x c_i'$ to
$T_i$ through product Lagrangian tori off $\Delta$ and the other five.
In $T^4$ the circles $b_1',b_1''$ are
disjoint and all four surgery circles avoid $\xi$. Thus the eight tori
are disjoint and miss the gluing region if the smoothing is small. After the surgeries in $\Sigma_2\x\Sigma_2$ and $T^4$, let
$\mathcal A$ and $\mathcal B$ be the respective exteriors of $\Delta$ and
$\Sigma'$. So $X=\mathcal A\cup_{\partial}\mathcal B$.

To pass from the spin example $X$ to its nonspin counterpart $X'$ while
preserving simple connectivity in the proofs of Propositions~\ref{p:groups} and~\ref{p:models} below, we perform one further $+1$ Luttinger
surgery, along the $\alpha$ factor of the torus $S=\alpha\x\beta$ in
Figure~\ref{fig:products}.
With the orientations shown, $[\alpha]=a_1+a_2$ and $[\beta]=d_1-d_2$.
For the ninth surgery we check $S$ against $\Delta$ and the six tori
in $P\x Q$, as the other two lie in $T^4$. Under $P\cong Q$, $\alpha$
runs east along $y=1$ in both handles and $\beta$ along $x=1$, north
in the first and south in the second. Near $C$ the arcs bend apart
(Figure~\ref{fig:tracks}, left). The collar map exchanges the horizontal
endpoints and fixes the vertical ones, so the smoothed circles are
disjoint and $S$ misses $\Delta$. The six product tori also miss $S$ throughout
the isotopy above, since their first circles stay in the annuli of
parallels disjoint from $\alpha$. The product symplectic form vanishes
on $S$, which therefore survives the sum and the eight surgeries as an
embedded Lagrangian torus of square zero in $X$. We can thus perform
the indicated Luttinger surgery on $S$ to obtain $X'$.
For this ninth torus we use a genus--$2$ geometric dual, obtained by
joining a punctured torus from each summand in the proof of
Lemma~\ref{lem:homology}.

Write $E=X\setminus\Int\nu S$, with $\mu$ the meridian of $S$ and
$\lambda$ the Lagrangian push-off of $\alpha$. Gluing $T^2\x D^2$
back to $E$ with its meridian identified with $\mu$ recovers $X$.
The Luttinger surgery giving $X'$ identifies this meridian with $\mu\lambda$.

\begin{lemma}\label{lem:homology}
Inclusion of the common exterior of the nine tori induces an
isomorphism onto $H_1(M)\cong\Z^8$.
Moreover, $H_1(E)=H_1(X)=H_1(X')=0$.
\end{lemma}

\begin{proof}
\linespread{1.10}\selectfont
It suffices to identify the eight surgery longitudes as a basis of $H_1(M)$
and show that the nine meridians vanish in $H_1$ of the common torus
exterior. To compute $H_1(M)$, observe that the slice $\{p\}\x Q$
meets the diagonal $\Delta$ once. On the $T^4$
side, choose $\zeta$ outside the projection of the smoothing region to
the first $T^2$. Then $\{\zeta\}\x T^2$ meets $\Sigma'$ once, at
$(\zeta,\xi)$. On either side the punctured dual is a surface in the
exterior with boundary the meridian. The meridian therefore vanishes
in $H_1$, and the exact sequence of the pair shows that inclusion of
the exterior induces an isomorphism on $H_1$. The two handles of $\Sigma'$ carry the four
coordinate generators of $H_1(T^4)$, so $H_1(\Sigma')\to H_1(T^4)$ is
an isomorphism. Mayer--Vietoris now gives $H_1(M)\cong H_1(P\x Q)\cong\Z^8$.

Write $a,b,c,d$ for the homology classes of $a_1,b_1,a_2,b_2$ in $P$ and
$p,q,r,s$ for those of $c_1,d_1,c_2,d_2$ in $Q$. For the based loop $d$ in
Section~\ref{sec:groups} we push the second point into $Q_2$. The surgery
longitudes give the basis
\begin{equation}\label{eq:basis}
                 p,\ q,\ r,\ s,\ a,\ c,\ b+q,\ d+s .
\end{equation}
To compute the meridians we choose parallels for the duals in~\eqref{eq:surgeries}.
For the six on the $P\x Q$ side, take the $b_i$
parallels at $x=6$ for $L_{12},L_{21}$, at $x=7$ for $L'_{12},L'_{21}$,
and at $x=8$ for $T_1,T_2$. The six dual tori are disjoint since their
projections to $P$ are disjoint. In $Q_j$ take the vertical $d_j$ at
$x=-6$ and the horizontal $c_j$ at $y=-6$.
The former crosses $c_j'$ once and misses $d_j'$, and the latter crosses
$d_j'$ once and misses $c_j'$. Thus each dual meets its designated
surgery torus once and misses the other seven. The dual of $T_i$ also
misses $\Delta$ since its two vertical circles have different $x$-coordinates.
For the two duals in $T^4$, choose disjoint horizontal circles in the
first $T^2$. In the second $T^2$ take a vertical parallel disjoint from
$b_2'$ for the dual of $T'_1$ and a horizontal parallel disjoint from
$a_2'$ for the dual of $T'_2$, with all four circles away from $\xi$
and the smoothing region. This gives eight mutually disjoint dual tori.
Only those dual to $L'_{12}$ and $L'_{21}$ meet
$S$, and each meets it once. The two intersection points have coordinates
$((7,1),(1,-6))$ in $P_1\x Q_2$ and $P_2\x Q_1$, respectively.

To obtain a dual of $S$, take $T_0$ in $P_1\x Q_1$ to be the product of the vertical $x=-7$ and
the horizontal $y=-3$. These circles cross $\alpha$ and $\beta$ once
each and meet once under $P_1\cong Q_1$. Thus $T_0$ meets $S$ and
$\Delta$ once each and misses the tori of~\eqref{eq:surgeries} and
their duals. On the $T^4$ side take $\{\zeta\}\x T^2$ with $\zeta$ off all
the chosen first circles and the smoothing disk. This torus meets
$\Sigma'$ once. Join the two punctured tori by a meridional annulus in
the neck. This gives a genus--$2$ surface $D_S$ disjoint from the
eight surgery tori and their duals. Its single transverse intersection with
$S$ is the point $((-7,1),(1,-3))$ on $T_0$, so $[S]$ is
primitive in $H_2(X)$.

The punctured $D_S$ gives $[\mu]=0$ in $H_1$ of the common exterior.
The twice-punctured duals to $L'_{12},L'_{21}$ give their meridians as
$\pm[\mu]$, and the other six duals have just one puncture. Thus all
nine meridians vanish in $H_1$. They generate the kernel of inclusion
on $H_1$ by the exact sequence of the pair, so inclusion is an
isomorphism. The eight surgeries then kill the basis~\eqref{eq:basis}
in $H_1(E)$, and hence in $H_1(X)$ and $H_1(X')$.
\end{proof}

Lemma~\ref{lem:homology} shows that the nine meridians vanish in
$H_1$ of the common exterior. In fact, they bound mutually disjoint
embedded oriented surfaces there. Puncturing the duals gives seven of
these at once. To obtain the other two, take two parallel copies of
$D_S$, using product push-offs of its two tori and a parallel annulus
in the neck. Pair each copy with one of the duals to $L'_{12},L'_{21}$,
orienting them to meet $S$ with opposite signs. Tube the two
intersection points together in $\partial\nu S$ along an arc in $S$,
taking the two arcs disjoint and away from $D_S\cap S$. Each resulting
genus--$3$ surface meets only its designated torus. Puncturing the dual
of genus $g$ thus expresses its meridian as a product of $g$ commutators
in the exterior group, up to conjugacy and orientation.

\linespread{.93}\selectfont
\begin{remark}\label{rem:model}
The model $M$ is \emph{nonaspherical}. Let $\Gamma\subset\pi_1(M)$ be
generated by the four boundary lifts on $\Sigma'$. Lifts from different
factors commute, and the Hopf band exterior identifies the other two
commutators with the central meridian $m_\Sigma$. The lifts span a rank
four subgroup of $H_1(M)=\Z^8$, so $\Gamma$ has infinite index and
$\Gamma/\gen{m_\Sigma}\cong\Z^4$. Were $M$ aspherical, $\Gamma$ would
be torsion-free and its noncompact cover would give $H_i(\Gamma)=0$
for $i\ge4$. Yet a $K(\Gamma,1)$ would be $T^4$ if $m_\Sigma=1$, or
the oriented circle bundle over $T^4$ given by the central extension
otherwise. Both have nonzero top homology.

Our isotopies put the six tori in the positions of~\cite[Addendum~12]{Fake},
whose decomposition shows that $\widehat{\mathcal A}=\mathcal A\cup\nu\Delta$
is aspherical and $C\x Q$ is $\pi_1$-injective there. On the gluing
boundary $[v,w]\ne1$, since killing $a_1,b_2$ in $\pi_1(\Sigma')$ leaves
$b_1,a_2$ free. On the $T^4$ side their product rectangle misses $\Sigma'$
and both surgery tori, so $[v,w]=1$ even after surgery. We use this
relation in Proposition~\ref{p:groups} to kill the remaining fiber generators.
\end{remark}

\section{The two exotic 4-manifolds}\label{sec:groups}

We now show that $X$ and $X'$ are indeed exotic copies of $\#_3(S^2\x S^2)$ and $\#_3(\CP\#\CPb)$, respectively.
We will first show that these manifolds are simply connected. From
$\pi_1(E)$, where $E=X\setminus\Int\nu S$, killing the meridian $\mu$
of $S$ gives $\pi_1(X)$ and killing $\mu\lambda$ gives $\pi_1(X')$.
Lemma~\ref{lem:meridian} shows that the diagonal meridian $t$ dies in each quotient.
Write $\mathcal A'$ for the diagonal exterior $\mathcal A$ after surgery
on $S$. Reattaching $\nu\Delta$ to either exterior adds just $t=1$ by
van Kampen, since the circle fiber $t$ generates the kernel of
$\pi_1(\partial\nu\Delta)\to\pi_1(\Delta)$. Thus
$\pi_1(\mathcal A)\to\pi_1(X)$ and $\pi_1(\mathcal A')\to\pi_1(X')$
factor through the respective groups with $\Delta$ restored.
We then show that $a,c$ from the two handles of $P$ are central and use
$[v,w]=1$ from Remark~\ref{rem:model} to kill the fiber generators
modulo $a,c$. Both groups are therefore abelian and hence trivial by
Lemma~\ref{lem:homology}.

We first describe $\pi_1(\mathcal B)$ on the $T^4$ side. Before
smoothing, the two coordinate tori have exterior the product of two
punctured tori. Cut the first along $b_1',b_1''$ into an annulus and a
pair of pants. Over both pieces the fiber is the second punctured
torus $V$, with a constant section off $a_2'\cup b_2'$. Across a cut the
complement is the cut circle times
$([-1,1]\x V)\setminus(\{0\}\x\gamma)$, where $\gamma=a_2'$ or $b_2'$.
Its halves retract to their end fibers and meet along the connected
$V\setminus\gamma$, so van Kampen adds no generators besides the end
fibers and the section. Carry the horizontal fiber generator across
$b_1'$, parallel to $a_2'$, and the vertical one across $b_1''$,
parallel to $b_2'$. The two paths differ by a loop in the section.
Thus the boundary lifts $u,v,w,z$ of $a_1,b_1,a_2,b_2$ on $\Sigma'$
generate, the first pair coming from the section and the second from
one fiber. Gluing back $T^2\x D^2$ in each Luttinger surgery adds only
relations, since $\pi_1(T^3)\to\pi_1(T^2\x D^2)$ is onto. In the smoothing ball the
Hopf band is a disk with one band, so its exterior is generated by one
meridian. Both positive coordinate meridians give this $m_\Sigma$,
central among the lifts as the boundary circle fiber.

The product square of coordinate loops $g,h$ in different factors has
counterclockwise boundary $[g,h]$. For $v,w$ and $v,z$ this is trivial
since the squares miss $\Sigma'$ and both surgery tori. At a transverse
intersection, take a small counterclockwise loop around the puncture,
based by an arc from the corner. The crossing is positive if the ordered
directions $g,h$ project to a positive basis of the normal plane oriented
before~\eqref{eq:surgeries}. The loop then gives the chosen meridian.
A negative basis gives the inverse meridian. At $T'_1$ the normals point west and north,
while $u,z$ run east and north, giving the inverse meridian. At $T'_2$
both normals point west, while $w,u^{-1}$ run east in the second factor
and west in the first. Exchanging the factors and reversing one direction
gives the positive meridian. The longitudes on the far side are
$uvu^{-1}$ and $wzw^{-1}$, respectively. Thus
\begin{equation}\label{eq:cap}
 [u,v]=[w,z]=m_\Sigma,\quad [v,w]=[v,z]=1,\qquad wu=uwz,\quad zu=uzv .
\end{equation}
The last two relations express $z$, then $v$, in terms of $u,w$.

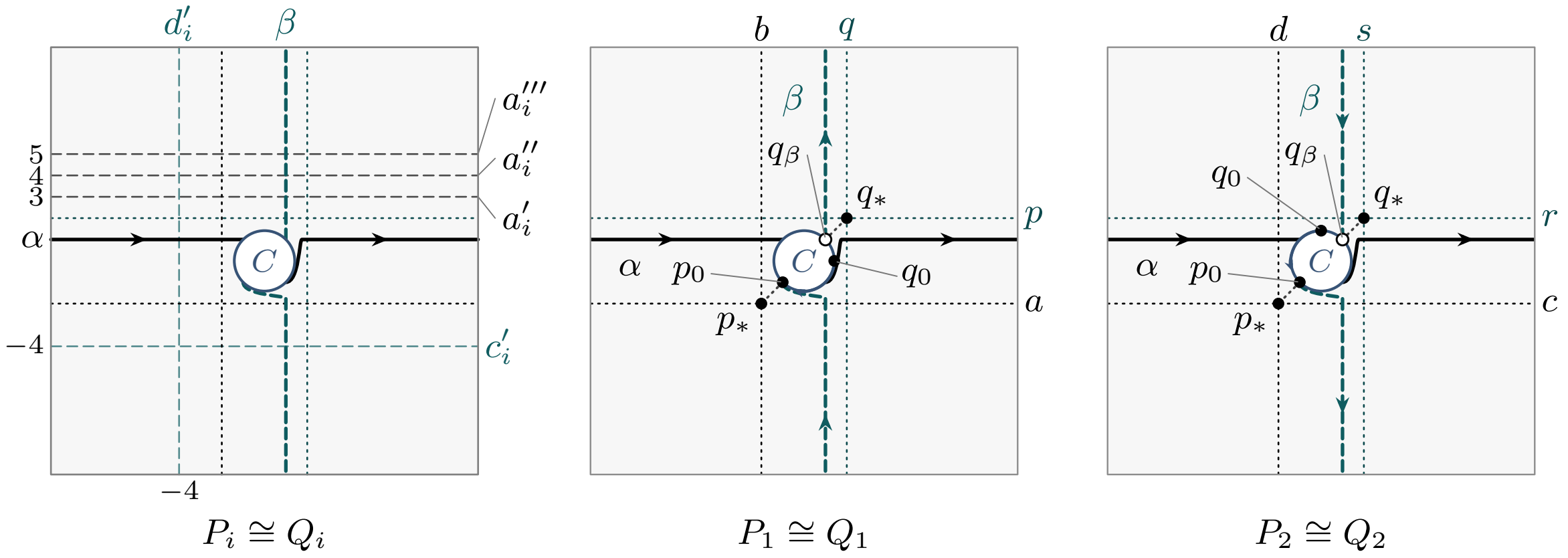
\begin{figure}[htb]
\colorlet{PInk}{black}
\definecolor{Qteal}{RGB}{15,92,96}
\definecolor{Cblue}{RGB}{55,83,118}
\tikzset{
  al/.style={draw=PInk, line width=0.9pt, line cap=round},
  be/.style={draw=Qteal, line width=0.9pt, dash pattern=on 2.6pt off 1.5pt},
  parP/.style={draw=PInk!70, line width=0.45pt, dash pattern=on 2.6pt off 1.5pt},
  parQ/.style={draw=Qteal!70, line width=0.45pt, dash pattern=on 2.6pt off 1.5pt},
  trP/.style={draw=PInk!80!black, line width=0.45pt, dash pattern=on 0.5pt off 1.4pt, line cap=round},
  trQ/.style={draw=Qteal!80!black, line width=0.45pt, dash pattern=on 0.5pt off 1.4pt, line cap=round},
  Cst/.style={draw=Cblue, line width=0.7pt, fill=white},
  rad/.style={draw=black!75, line width=0.6pt, dash pattern=on 0.4pt off 1.3pt, line cap=round},
  pt/.style={circle, fill=black, inner sep=1.05pt},
  opt/.style={circle, draw=black, fill=white, line width=0.5pt, inner sep=1.05pt},
  arr/.style={-{Stealth[length=3.05pt,width=2.2pt]}},
  point label/.style={draw=black!53,line width=.33pt},
  marr/.style={postaction={decorate}, decoration={markings, mark=at position #1 with {\arrow{Stealth[length=4.2pt,width=3.4pt]}}}},
  marrr/.style={postaction={decorate}, decoration={markings, mark=at position #1 with {\arrowreversed{Stealth[length=4.2pt,width=3.4pt]}}}},
  lab/.style={font=\footnotesize},
}
\def\seamradius{1.41421356}
\newcommand{\alphaturn}{(1,-1) .. controls (1.625,-1) and (1.625,1) .. (1.75,1)}
\newcommand{\betaturn}{(-1,-1) .. controls (-1,-1.625) and (1,-1.625) .. (1,-1.75)}
\newcommand{\handlepaths}[5]{
 \begin{scope}
  \clip (-10,-10) rectangle (10,10);
  \fill[black!3] (-10,-10) rectangle (10,10);
  \draw[trP] (-2,-10)--(-2,10); \draw[trP] (-10,-2)--(10,-2);
  \draw[trQ] (2,-10)--(2,10);   \draw[trQ] (-10,2)--(10,2);
  \draw[rad] (-1,-1)--(-2,-2); \draw[rad] (1,1)--(2,2);
  \draw[al,marr=.44] (-10,1)--(-1,1);
  \draw[al,marr=.73] \alphaturn -- (10,1);
  \ifnum#1=1
    \draw[be,marr=.57] (1,1)--(1,10);
    \draw[be,marrr=.73] \betaturn -- (1,-10);
  \else
    \draw[be,marrr=.57] (1,1)--(1,10);
    \draw[be,marr=.73] \betaturn -- (1,-10);
  \fi
  \draw[Cst] (0,0) circle (\seamradius);
  \ifnum#1=1
    \draw[Cblue,line width=.75pt,arr] (0:\seamradius) arc (0:-112:\seamradius);
    \node[pt] at (0:\seamradius) {};
  \else
    \draw[Cblue,line width=.75pt,arr] (90:\seamradius) arc (90:202:\seamradius);
    \node[pt] at (90:\seamradius) {};
  \fi
  \node[pt] at (-1,-1) {}; \node[opt] at (1,1) {};
  \node[pt] at (-2,-2) {}; \node[pt] at (2,2) {};
 \end{scope}
 \draw[black!45,line width=.4pt] (-10,-10) rectangle (10,10);
 \node[font=\scriptsize,text=Cblue,inner sep=0pt] at (0,0) {$C$};
 \draw[point label] (-1.10,-.98) -- (-4.35,-.60);
 \node[lab,anchor=east,inner sep=.8pt] at (-4.35,-.60) {$p_0$};
 \draw[point label] (.98,1.12) -- (.15,4.95);
 \node[lab,anchor=east,inner sep=.9pt] at (.15,4.95) {$q_\beta$};
 \ifnum#1=1
   \draw[point label] (1.55,-.05) -- (4.35,-.72);
   \node[lab,anchor=west,inner sep=.8pt] at (4.35,-.72) {$q_0$};
 \else
   \draw[point label] (-.10,1.52) -- (-3.45,3.90);
   \node[lab,anchor=east,inner sep=.8pt] at (-3.45,3.90) {$q_0$};
 \fi
 \node[lab,anchor=north east,inner sep=1.65pt] at (-2.10,-2.10) {$p_*$};
 \node[lab,anchor=south west,inner sep=1.65pt] at (2.10,2.10) {$q_*$};
 \node[lab,text=PInk,anchor=north] at (-8.15,.75) {$\alpha$};
 \node[lab,text=Qteal,anchor=east] at (.72,7.50) {$\beta$};
 \node[lab,text=PInk!80!black,anchor=west,inner sep=1.5pt] at (10,-2) {$#2$};
 \node[lab,text=PInk!80!black,anchor=south,inner sep=1.5pt] at (-2,10) {$#3$};
 \node[lab,text=Qteal!80!black,anchor=west,inner sep=1.5pt] at (10,2) {$#4$};
 \node[lab,text=Qteal!80!black,anchor=south,inner sep=1.5pt] at (2,10) {$#5$};
}
\centering
\resizebox{\linewidth}{!}{\begin{tikzpicture}[x=.19cm,y=.19cm,line cap=round,line join=round]
 \begin{scope}[shift={(.38cm,.38cm)}]
  \fill[black!3] (-10,-10) rectangle (10,10);
  \foreach \y/\ylab/\l in {3/2.0/{a_i'},4/4.8/{a_i''},5/7.6/{a_i'''}}{
    \draw[parP] (-10,\y)--(10,\y);
    \draw[PInk!55,line width=.3pt] (10,\y)--(10.8,\ylab);
    \node[lab, text=PInk, anchor=west, inner sep=1.5pt] at (10.8,\ylab) {$\l$};
    \node[font=\tiny, anchor=east, inner sep=1.5pt] at (-10,\y) {$\y$}; }
  \draw[parQ] (-10,-4)--(10,-4); \node[lab, text=Qteal, anchor=west, inner sep=1.5pt] at (10,-4) {$c_i'$};
  \draw[parQ] (-4,-10)--(-4,10); \node[lab, text=Qteal, anchor=south, inner sep=1.5pt] at (-4,10) {$d_i'$};
  \node[font=\tiny, anchor=east, inner sep=1.5pt] at (-10,-4) {$-4$};
  \node[font=\tiny, anchor=north, inner sep=1.5pt] at (-4,-10) {$-4$};
  \draw[trP] (-2,-10)--(-2,10); \draw[trP] (-10,-2)--(10,-2);
  \draw[trQ] (2,-10)--(2,10);   \draw[trQ] (-10,2)--(10,2);
  \draw[al, marr=0.5] (-10,1)--(-1,1); \draw[al, marr=0.6] \alphaturn -- (10,1);
  \draw[be] (1,1)--(1,10); \draw[be] \betaturn -- (1,-10);
  \draw[Cst] (0,0) circle (\seamradius);
  \node[font=\scriptsize,text=Cblue,inner sep=0pt] at (0,0) {$C$};
  \draw[black!50, line width=0.45pt] (-10,-10) rectangle (10,10);
  \node[lab, text=PInk, anchor=east, inner sep=1.5pt] at (-10,1) {$\alpha$};
  \node[lab, text=Qteal, anchor=south, inner sep=1.5pt] at (1,10) {$\beta$};
  \node[lab, anchor=north] at (0,-11.4) {$P_i\cong Q_i$};
 \end{scope}
 \begin{scope}[shift={(5.18cm,.38cm)}]
  \handlepaths{1}{a}{b}{p}{q}
  \node[lab,anchor=north] at (0,-11.4) {$P_1\cong Q_1$};
 \end{scope}
 \begin{scope}[shift={(9.78cm,.38cm)}]
  \handlepaths{2}{c}{d}{r}{s}
  \node[lab,anchor=north] at (0,-11.4) {$P_2\cong Q_2$};
 \end{scope}
\end{tikzpicture}}
\caption{On the left we draw the parallel surgery curves in the standard
square model $[-10,10]^2$ for $P_i\cong Q_i$. The middle and right squares
give the basing paths for $i=1,2$. Dotted radii join $p_0$ to $p_*=(-2,-2)$
and $q_\beta$ to $q_*=(2,2)$, the base points of the dotted coordinate
circles. The arrow on $C$ runs from $q_0$ toward $p_0$.}
\label{fig:tracks}
\end{figure}

We base the loops as in Figure~\ref{fig:tracks}. The ends of $\alpha$
on $C$ are $(-1,1),(1,-1)$ and those of $\beta$ are
$q_\beta=(1,1),p_0=(-1,-1)$. Starting at $o=(p_0,q_0)$, with
$q_0$ the rightmost point of $C$ in the coordinates of $Q_1$, move the second point along $C$
to $q_\beta$, then push both radially to $p_*,q_*$. For $a,b,p,q$ in
$P_1\x Q_1$ and $c,d,r,s$ in $P_2\x Q_2$, take the indicated
coordinate circle with the other point fixed and return by the same
path. For $d$ in $P_2$, write $\bar d$ for the loop with the second
point left at $q_0$ or pushed into $Q_1$.

To base the diagonal meridian $t$, follow $C$ from $q_0$ past
$(1,-1)$ to just before $p_0$ and turn once positively about $p_0$.
For the lifts $A,B,W,D$, bring the second point close to $p_0$ by
the same arc and push the pair radially until the first reaches $p_*$.
Move the second along the shorter arc about $p_*$ to lie due east
in $Q_1$ or due north in $Q_2$. Keep this displacement as the
first point goes around $a_1,b_1,a_2,b_2$, respectively, and return
by the basing path.

We choose the gluing over the identification $\Sigma'\cong\Delta$
of Section~\ref{sec:construction} so that
\begin{equation}\label{eq:flat}
        m_\Sigma=t^{-1},\qquad (u,v,w,z)=(tA,t^{-1}B,W,t^{-1}D).
\end{equation}
Indeed, the opposite normal Euler numbers $2,-2$ give a circle bundle
identification reversing the fiber, with $[u,v][w,z]=m_\Sigma^2$ and
$[A,B][W,D]=t^{-2}$. Lifts of the same loop differ by a power of the
central meridian. Rotating the fibers by a map $\Delta\to S^1$ sets
these powers independently on the four basis loops, as in~\eqref{eq:flat}.

\begin{lemma}\label{lem:boundary-lifts}
In $\pi_1(E,o)$, $a$ commutes with $p,r,s$ and $c$ with $p,q,r$.
With $\mu,\lambda$ based by the first rectangle in Figure~\ref{fig:rectangles},
the push-off of $\alpha$ is $\lambda=ac$ and the relations~\eqref{eq:cross} hold.
For $k=sr$ the boundary lifts satisfy
\[
                 A=(ap)^q,\qquad w=(cr)^k,\qquad D=(ds)^r .
\]
Moreover, $a^q=t^{-1}a$ and $c^k=ct^{-1}$.
\end{lemma}

\begin{proof}
We read the commutations, the two peripheral curves of $S$, and finally
the diagonal lifts from the based paths in Figure~\ref{fig:tracks}.
For each commuting pair, the horizontal $a$ or $c$ at $y=-2$ and
its basing radius miss $\alpha$ and the parallels at $y=3,4,5$, so
the product square misses $S$ and the six surgery tori. It misses
$\Delta$ as well. The two circles lie in different handles or on the
distinct horizontals $y=-2,2$, and the paths $p_0p_*$ and
$q_0q_\beta q_*$ miss each other and the opposite circles.

For $S$, run the second point west from $q_*$ and the first north
from $p_*$. This gives the two rectangles in Figure~\ref{fig:rectangles},
with positions $((-2,-2+\tau),(2-\rho,2))$ in $P_2\x Q_i$ for
$0\le\rho,\tau\le20$. Base the small loops by straight segments from
$(0,0)$ to the punctures. The counterclockwise boundary reads the lower
one first, then the one to $S$ at $(1,3)$.

\begin{figure}[htb]
\centering
\begin{tikzpicture}[x=1cm,y=1cm,line cap=round,line join=round,
                    font=\footnotesize,inner sep=1.3pt]
  \definecolor{tailInk}{RGB}{55,83,118}
  \tikzset{forward/.style={postaction={decorate},decoration={markings,
      mark=at position #1 with {\arrow{Stealth[length=1.40mm,width=1.02mm]}}}}}
  \def\loopradius{2.8mm}
  \def\puncture#1#2#3#4{
    \coordinate (hole) at ({.70*#1},{.21*#2});
    \draw[tailInk,line width=.6pt,shorten >=\loopradius] (0,0)--(hole);
    \draw[tailInk,line width=.9pt] (hole) circle[radius=\loopradius];
    \coordinate (east) at ($(hole)+(0:\loopradius)$);
    \fill[tailInk] ($(east)+(0,1.25mm)$) -- ($(east)+(-1.0mm,-.85mm)$)
      -- ($(east)+(1.0mm,-.85mm)$) -- cycle;
    \node[text=tailInk,anchor=west,inner sep=0pt]
      at ($(east)+(1.8mm,.65mm)$) {$#4$};
    \fill (hole) circle[radius=.76pt];
    \node[above=1.2pt] at ($(hole)+(90:\loopradius)$) {$#3\;{\scriptstyle (#1,#2)}$};}
  \foreach \offset/\handle in {0/2,7.27/1}{
    \begin{scope}[shift={(\offset,0)}]
      \path[fill=black!2] (0,0) rectangle (5.45,1.98);
      \draw[forward=.82] (0,0)--(5.45,0);
      \draw[forward=.59] (5.45,0)--(5.45,1.98);
      \draw[forward=.51] (5.45,1.98)--(0,1.98);
      \draw[forward=.55] (0,1.98)--(0,0);
      \ifnum\handle=2
        \def\baseEdge{r^{-1}}\def\otherEdge{d}
        \puncture{1}{3}{S}{\mu}
        \puncture{4}{4}{\Delta}{t}
      \else
        \def\baseEdge{p^{-1}}\def\otherEdge{\bar d}
        \puncture{1}{3}{S}{\mu^{-1}}
        \puncture{6}{6}{L'_{21}}{\mu_{L'_{21}}^{-1}}
      \fi
      \fill (0,0) circle[radius=1.42pt] node[below left=1pt] {$(0,0)$};
      \node[below=1.65pt] at (2.71,0) {$\baseEdge$};
      \node[right=2pt] at (5.45,.99) {$\otherEdge$};
      \node[above=2pt] at (2.72,1.98) {$P_2\x Q_{\handle}$};
    \end{scope}}
\end{tikzpicture}
\begin{equation}\label{eq:rectangles}
 \makebox[.46\linewidth]{\([r^{-1},d]=t\mu,\)}
 \makebox[.46\linewidth]{\([p^{-1},\bar d]=(\bar d q^{-1}\bar d^{-1})^p\,\mu^{-1}.\)}
\end{equation}
\caption{The rectangles giving~\eqref{eq:rectangles}, with outer parts
shortened. All four small loops run counterclockwise. At $S=\alpha\x\beta$,
$\beta$ runs south in $Q_2$ and north in $Q_1$, so its chosen normal
points east in $Q_2$ and west in $Q_1$. Thus the loop about $S$ gives
$\mu$ in the first rectangle and $\mu^{-1}$ in the second. The other
loops give the meridian $t$ of $\Delta$ and the inverse meridian
$\mu_{L'_{21}}^{-1}$ of $L'_{21}$.}
\label{fig:rectangles}
\end{figure}

In $P_2\x Q_2$ the other puncture is at $(4,4)$, where the two
points coincide at $(-2,2)$. The normal coordinates $(4-\rho,4-\tau)$
at $\Delta$ and $(\tau-3,1-\rho)$ at $S$ both preserve the orientation
of the $(\rho,\tau)$ rectangle. In $P_2\x Q_1$ the reversal of $\beta$
changes $1-\rho$ to $\rho-1$ at $S$, reversing the normal orientation.
At $L'_{21}=a_2''\x d_1'$ the normal coordinates $(\tau-6,\rho-6)$
also reverse orientation by interchanging the parameters. The radius from $q_\beta$ leaves on the left
of $\beta$ in $Q_2$ and on the right in $Q_1$, giving $\bar d=\mu d$.

The segment to $(1,3)$ in the first rectangle bases the push-off
$\lambda$ just east of $\beta$ in $Q_2$. From $(-2,1)$ in $P_2$ it
follows $\alpha$ through $P_1$ and back through $P_2$, so $\lambda=ac$.
In the second rectangle push the negative $d_1'$ longitude above and
to the right of $(6,6)$ to the far corner. The bottom and right edges
reach it by $p^{-1}\bar d$, giving the conjugate of $q^{-1}$
in~\eqref{eq:rectangles}. The negative $b_2,d_2$ rectangle meets only
$T_2$, giving $[d^{-1},s^{-1}]=c^{-1}$. With $\bar d=\mu d$, the
three rectangle relations now read
\begin{equation}\label{eq:cross}
                  d^r=t\mu d,\qquad pd=\mu dqp,\qquad s^d=sc .
\end{equation}

To split the simultaneous motions $A,w=W,D$ into coordinate loops,
move the second point south, west then south, and west, respectively.
The basing paths read $q^{-1},r^{-1}s^{-1}=k^{-1},r^{-1}$, giving the
conjugates $(ap)^q,(cr)^k,(ds)^r$ on return. The rectangles
$[q^{-1},a]$ and $[c^{-1},k^{-1}]$ each meet $\Delta$ once negatively.
For the first, move the second point south and the first east. For the
second, move the second point west to the first vertical, both west,
then the second south. Sliding the interval between the points back
along either path identifies the positive meridian with $t$. Thus both
commutators equal $t^{-1}$, giving $a^q=t^{-1}a$ and $c^k=ct^{-1}$.
\end{proof}

Read the punctured-torus boundaries in the second fiber from the $q_0$
side. In order they give $t=\delta_2\delta_1$, where
$\delta_1=p^{-1}q^{-1}pq$ and $\delta_2=r^{-1}s^{-1}rs$.

\begin{lemma}\label{lem:meridian}
The diagonal meridian $t$ is trivial in both $\pi_1(X)$ and $\pi_1(X')$.
\end{lemma}

\begin{proof}
Since $a^q=t^{-1}a$ by Lemma~\ref{lem:boundary-lifts}, our aim is to
show $[a,q]=1$. We first work in $\pi_1(E)$, where we retain the meridian $\mu$
of $S$. With $h=s^r$, $y=r^k=\delta_2r$ we have $k=rh$, $\delta_2=h^{-1}s$, and the lemma and~\eqref{eq:flat} give
\[
\begin{alignedat}{2}
  u &=t(ap)^q=ap\delta_1, &\qquad A &=ut^{-1}=ap\delta_2^{-1},\\
  tw&=yc,                 & t &=h^{-1}h^c .
\end{alignedat}
\]
For the second line use $w=r^kc^k=yc t^{-1}$ and $t=k^{-1}k^c$.

To compute $f(c)$ for $f(g)=DgD^{-1}$, observe that $D=r^{-1}dk=d^rh$
and $Dc=(sd)^r=hd^r$ give $f(ch)=h$, while~\eqref{eq:cross} gives
$f(y)=t\mu r$. Applying $f$ to $w=ct^{-1}y$ and using $f(t)=t$,
$f(w)=tw$ gives $f(c)\mu r=tw=yc$, whence
\begin{equation}\label{eq:transport}
                 f(c)=\delta_2c\mu^{-1},\qquad Dc=ht\mu d .
\end{equation}

In the off-diagonal blocks the $q,s$ longitudes give
$[q,\bar d]=p\mu^{-1}$ and $[s,b]=s\mu^as^{-1}r$. In $X$ the meridian $\mu$ bounds
a normal disk to $S$, so $\mu=1$ and $\bar d=d$. We recover the two twists
of~\cite[Theorem~9.2 and Section~9.2.1]{Kirk} in $P_2\x Q_1$ with
factors exchanged and second fiber longitude $zyz^{-1}$~\cite[(4)]{Fake}.
Indeed, $p^d=qp$ by~\eqref{eq:cross} and $qd=pdq$ gives $q^d=qpq$.
Exchanging handles~\cite[Lemma~4]{Fake} gives
\[
         (p^d,q^d)=(qp,qpq),\qquad (r^b,s^b)=(sr,srs).
\]
In particular $d$ fixes $\delta_1$ and $dpd^{-1}=q^{-1}p^2$.
In the relation $Ay\,cD=yc\,u$ from~\eqref{eq:cap}, substitute
$cD=\delta_2^{-1}h\delta_2\,d\delta_1$ from~\eqref{eq:transport}
and $u=ap\delta_1$. Cancelling $\delta_1$ gives $d=j\,c\,a\,p$ with
$j=(Ay\delta_2^{-1}h\delta_2)^{-1}y\in\gen{A,r,s}$. This subgroup
centralizes $a$, and $c\,a\,p$ centralizes $p$. Thus $dpd^{-1}=jpj^{-1}$
commutes with $a$. Since $dpd^{-1}=q^{-1}p^2$ and $[a,p]=1$, we get
$[a,q]=1$.

In $X'$ we have $\mu=(ac)^{-1}$. Our aim is to put $q$ in $\gen{a,H}$,
where $H=\gen{p,r,s}=\gen{p,r,h}$ is centralized by $a$. To see that
$c$ normalizes $H$, it suffices to check the two conjugates of $h$, since
$c$ fixes $p,r$. Applying $f$ to $t=(ch)^{-1}hc$ and using $f(t)=t$, $f(ch)=h$
gives
\[
             c^{-1}hc=ht=f(c)^{-1}h f(c).
\]
It follows that $f(c)c^{-1}=h^{-1}sca$ commutes with $h$, and hence so
does $sc$, since $a$ centralizes $H$. This gives $chc^{-1}=s^{-1}hs\in H$.
For the other conjugate, $yc=tw$ commutes with $t$, so
$t=y(ctc^{-1})y^{-1}=y(chc^{-1})^{-1}hy^{-1}\in H$ and $c^{-1}hc=ht\in H$.

Now $u=tap\delta_2^{-1}$ and $w=t^{-1}yc$ normalize $H$, as do $D=tz$
and $d=rDk^{-1}$ since $z\in\gen{u,w}$ by~\eqref{eq:cap} and $t,r,k\in H$.
Thus $c=s^{-1}s^d$ belongs to $H$. The same expressions put $u,w,D,d$ in
$\gen{a,H}$, and then $pd=\mu dqp$ with $\mu=(ac)^{-1}$ puts $q$ there as well.
This proves $[a,q]=1$.
\end{proof}

\begin{prop}\label{p:groups}
The manifolds $X$ and $X'$ are simply connected.
\end{prop}

\begin{proof}
Let $G=\pi_1(X)$ or $\pi_1(X')$. With $\Delta$ restored on the $P\x Q$
side as above, we show that $a,c$ are central and generate $G$. Cut $P$ along the first circles of the surgery tori.
Over each piece $U$ the diagonal and one fiber generate $\pi_1(U\x Q)$.
Each gluing face contains a whole fiber and the gluing fixes the diagonal,
since the two circles of each torus are disjoint under $P\cong Q$.
The remaining van Kampen generators lift across the collars as
$(\gamma,\gamma)$, so the diagonal and one fiber generate $G$~\cite[Section~2]{Fake}.

Since the diagonal loops come from $\pi_1(\mathcal B)=\gen{u,w}$,
where $u=ap\delta_1$ and $w=(cr)^k$, we have $G=\gen{a,c,F_Q}$ with
$F_Q=\gen{p,q,r,s}$. With $t=1$, $a^q=a$ and $c^{sr}=c$
give the two commutations missing from Lemma~\ref{lem:boundary-lifts}.
Thus $F_Q$ is normal and $c=s^{-1}s^d\in F_Q$, so $a,c$ are central.

It remains to kill the fiber in $N=G/\gen{a,c}$. Here $\mu=1$, so
the twists for $X$ apply and~\eqref{eq:transport} gives $\delta_2=1$,
hence $\delta_1=1$. Thus $[p,q]=[r,s]=1$, $u=p$, $w=r$ and $z=D=ds$,
where $d$ commutes with $r,s$ by~\eqref{eq:cross}. The relation
$wu=uwz$ becomes $r^p=rz$. Conjugating by $d$, which fixes $r,z$ and
sends $p$ to $qp$, we get $r^{qp}=r^p$, so $[q,r]=1$.

To read $v=B$ in $N$, cap the boundaries of $P_1,Q_1$ using $a=1$
and $[p,q]=1$. In the $b,q$ square we reach $a_1'''$ before $c_1'$,
so $T_1$ meets the upper triangle and the lower reads $bq$ in $P_1\x Q_1$.
Basing conjugates $bq$ by $q^{-1}pq=p$ in $N$. The $b,p$ square misses
$T_1$ since its horizontal side misses $c_1'$. Its meridians at $\Delta,S$
die with $t,\mu$, giving $[b,p]=1$. As $[p,q]=1$ too, we have $v=bq$ even
with $T_1$ removed. Now $r=w$ commutes with $v$ by~\eqref{eq:cap} and with $q$, hence with
$b$. The twists $r^b=sr$, $s^b=srs$ then give $s=r=1$, so \eqref{eq:cap}
gives $z=v=1$. With $d=zs^{-1}=1$, the other two twists $p^d=qp$,
$q^d=qpq$ give $q=p=1$. Thus $N=1$ and $G$ is abelian on $a,c$,
so Lemma~\ref{lem:homology} gives $G=1$.
\end{proof}

Torus surgery preserves $\eu(M)=\eu(P\x Q)-2\eu(\Delta)=8$ and
$\sigma(M)=0$. Thus $b_2^+=b_2^-=3$, and it remains to determine parity.

\begin{prop}\label{p:models}
The manifolds $X$ and $X'$ are irreducible symplectic 4-manifolds
homeomorphic to $\#_3(S^2\x S^2)$ and $\#_3(\CP\#\CPb)$, respectively.
\end{prop}

\begin{proof}
We construct a spin structure on $X$ and then show that its restriction
to $E$ does not extend over the surgery on $S$. Any spin structure on $X'$
would restrict to this one, since $H^1(E;\Z_2)=0$ by Lemma~\ref{lem:homology}.
The product spin structures agree on the normal fiber and differ by a
class from $H^1(\Delta;\Z_2)$. Since $H^1(T^4;\Z_2)\to H^1(\Sigma';\Z_2)$
is onto, they can be matched, so $M$ is spin.

In the Lagrangian framing let $\eta$ be the difference from the invariant boundary
spin structure. A primitive straight slope bounds a spin disk exactly when
its $\eta$-value is $1$. By~\eqref{eq:basis} a class in $H^1(M;\Z_2)$ sets
$\eta=0$ on the eight longitudes, so $\eta=1$ on the new meridians and $X$ is spin.

For the surgery on $S$, let $q$ be the quadratic refinement on
$H_1(P;\Z_2)$, with value one on a nonbounding spin circle.
The second-factor framing is constant along a first longitude, so
$\eta(\gamma')=1-q([\gamma])$. Thus $\eta=0$ on $a,c$ gives
$q(a_1)=q(a_2)=1$. Since $a_1\cdot a_2=0$,
$q([\alpha])=q(a_1+a_2)=0$. Hence
$\eta(\lambda)=1$ and $\eta(\mu+\lambda)=0$, making $X'$ nonspin.

Hence the homeomorphism types are as claimed by Freedman~\cite{Freedman}.
The spin manifold $X$ is minimal, and the Luttinger surgery giving
$X'$ preserves minimality~\cite{HoLi,Welschinger}.
Both are irreducible by~\cite{HK}
and hence exotic copies of the reducible smooth 4-manifolds
$\#_3(S^2\x S^2)$ and $\#_3(\CP\#\CPb)$, respectively.
\end{proof}

\linespread{1.08}\selectfont
\section{Proofs of Theorem~\ref*{thm:main} and its corollaries}\label{sec:families}

We first prove Theorem~\ref{thm:main} by producing infinitely many
pairwise nondiffeomorphic irreducible smooth structures in each of the
two homeomorphism types of Proposition~\ref{p:models}, infinitely many
symplectic and infinitely many nonsymplectic. To do this by knot surgery~\cite{FS},
we need a symplectic torus with simply connected exterior in each example, so once again, the $\pi_1$ computation is the most delicate issue to address here.

Perform the seven surgeries
in~\eqref{eq:surgeries} away from $T_1$, and let $E_1$ be its exterior
at this stage. Next do the Luttinger surgery on $S$, giving a spin
manifold $Y_0$ in which $T_1$ has exterior $E'_1$. The final surgery
on $T_1$ gives $Y_1=X'$, where the surgery core has the same exterior $E'_1$.

\begin{lemma}\label{lem:exteriors}
The group $\pi_1(E_1)$ is infinite cyclic on $a$, and $E'_1$ is simply connected.
\end{lemma}

\begin{proof}
We repeat the argument of Proposition~\ref{p:groups} with $T_1$
removed, then compute the abelianizations. The rectangles used to
kill $t$ miss $T_1$. After putting back the diagonal, cut $P$ along
$a_1'''$. The cut surface is connected, and the two end fibers generate
across the collar since $Q\setminus c_1'$ is connected, just as at
the start of Section~\ref{sec:groups}. The cut circle and a loop
crossing it once both lift to the diagonal, so the diagonal and one
fiber still generate. The lower triangle giving $v=bq$ and the $b,p$
square both miss $T_1$. The proposition therefore still makes $a,c$
central and kills the fiber modulo $\gen{a,c}$.

Thus both groups are abelian, and all meridians vanish in homology by
Lemma~\ref{lem:homology}. The seven longitudes in~\eqref{eq:basis}
leave only $a$, while the additional surgery on $S$ kills $a+c=a$.
\end{proof}

\begin{proof}[Proof of Theorem~\ref{thm:main}]
Take $Y=Y_0$ or $Y_1$, and let $T$ be $T_1$ or its surgery core,
respectively. In either case the exterior $E'_1$ is simply connected
by Lemma~\ref{lem:exteriors}, so $Y$ is simply connected and $T$ has
an integral dual. The spin calculation of Proposition~\ref{p:models},
with $a+c$ in place of $a$ in~\eqref{eq:basis}, gives a minimal spin
copy $Y_0$ of $\#_3(S^2\x S^2)$, while $Y_1=X'$ is our minimal
nonspin example. Perturb the form to make the homologically essential
Lagrangian torus $T$ symplectic~\cite{Gompf}. In the knot surgery the
simply connected exterior kills both the $S^1$ factor of
$S^1\x(S^3\setminus\Int\nu K)$ and a meridian of $K$, which normally
generates the knot group. Thus $Y_K$ is simply connected, and since knot
surgery preserves the intersection form, $Y_K$ is homeomorphic to $Y$
by Freedman's theorem~\cite{FS,Freedman}. Since $b_2^+(Y)=3$, the
Fintushel--Stern formula~\cite{FS,FSlectures} gives
$\mathrm{SW}_{Y_K}=\mathrm{SW}_Y\,\Delta_K(\exp(2\,\mathrm{PD}[T]))$
in the integral group ring of $H^2(Y;\Z)$, where $\Delta_K$ is the
symmetrized Alexander polynomial with $\Delta_K(1)=1$.

For the symplectic family take $K=\#_nJ$ with $J$ fibered and
$\Delta_J\ne1$. The knot surgeries are minimal symplectic sums, hence
irreducible~\cite{Usher,HK}. Since $\mathrm{SW}_Y\ne0$ and
$\Delta_K=\Delta_J^n$, the number of irreducible factors of
$\mathrm{SW}_{Y_K}$, counted with multiplicity, increases with $n$.
A fixed smooth manifold has at most two values of this count, one for
each orientation, so we obtain infinitely many diffeomorphism types.

For the nonsymplectic family take the twist knots $K_j$ with
$\Delta_{K_j}(s)=1+j(2-s-s^{-1})$ and $j\ge2$ even.
These manifolds are irreducible by the standard knot surgery argument
(see, e.g.,~\cite{AHP}). Fix a real cohomology
class pairing nontrivially with $[T]$. All terms of largest pairing in
the knot surgery formula come from the same end of $\Delta_{K_j}$,
so their coefficients are $-j$ times the corresponding coefficients
of $\mathrm{SW}_Y$. Thus the largest absolute invariant is unbounded
with $j$, distinguishing an infinite subfamily. If $Y_{K_j}$ admitted a symplectic
form, we could perturb it to pair nontrivially with $[T]$. The same
divisibility holds in the orientation induced by the form, since
$b_2^+(Y)=b_2^-(Y)=3$ and mirroring preserves the Alexander polynomial.
Yet the canonical class maximizes this pairing and has coefficient
$\pm1$ by Taubes~\cite{Taubes,TaubesGr}.

\end{proof}

\begingroup
\linespread{1.04}\selectfont
\smallskip
We now generalize our main result to any odd number of summands greater than one.

\begin{corollary}\label{cor:larger}
For all odd $m,n>3$, both $\#_m(S^2\x S^2)$ and $\#_n(\CP\#\CPb)$
admit infinitely many pairwise nondiffeomorphic irreducible smooth
structures, infinitely many of them supporting symplectic structures,
and infinitely many supporting none.
\end{corollary}

\begin{proof}
To increase $b_2^\pm$ by $2g$ while keeping the knot surgery torus $T$
with simply connected exterior, we sum along a genus--$2$ surface $F$
disjoint from $T$. Interchanging the factors of $T_0$ gives a torus
$A_0$ missing the six product tori and $S$, since its second circle
misses $\beta$. Take $B_0=\{\zeta\}\x T^2$ just right of $\xi$, off
the smoothing disk and surgery circles. To make $A_0$ symplectic,
choose closed $1$-forms $f_i,g_i$ on $P,Q$ supported in handle $i$,
with periods $(1,0)$ on $(a_i,b_i)$ and $(0,1)$ on $(c_i,d_i)$,
respectively. Then $-\sum_i\operatorname{pr}_1^*f_i\wedge
\operatorname{pr}_2^*g_i$ integrates negatively on $A_0$ and to zero
on all seven product tori, since $[\alpha]=a_1+a_2$ and
$[\beta]=d_1-d_2$. An exact adjustment makes it vanish near those tori
and restrict to a negative area form on $A_0$. A small perturbation
makes $A_0$ symplectic with positive intersection with $\Delta$, so
it joins $B_0$ to give the symplectic $F$ of square zero. We make
$T$ symplectic using its dual.

We need one relation on the generators $w_0,z_0$ of $B_0^\circ\subset F$.
Push them towards $\xi$ to a parallel fiber. Before smoothing, over
a disk containing $\xi$ and the parallel point, off $\zeta$ and the
surgeries, the exterior is an annulus times a punctured torus. The
annulus circle commutes with $w_0$ and becomes the meridian
$[w_0,z_0]^{\pm1}$ on smoothing, giving $[[w_0,z_0],w_0]=1$ in the
exterior of $F\cup T$.

Over $B_0^\circ$ take the $\Sigma_g$-bundle, $g\ge1$, with
monodromies $U=\prod_i t_{x_i}$, $V=\prod_i t_{y_i}$ for the standard
pairs $x_i,y_i$, using area-preserving twists fixing a disk. On the
pullback $V_g$ by the fold $F\to B_0^\circ$, the invariant fiber area
form plus an area form from $F$ is symplectic. A point of the fixed disk
gives a square zero symplectic section $S_g$. Form $Z_g=Y\mathbin{\#}_{F=S_g}V_g$,
matching the parallels of $w_0,z_0$ to those from the fixed disk.
The relation just found makes $J=[U,V]$ commute with $U$ on the image
of the punctured fiber in $\pi_1(Z_g\setminus\nu T)$. On each pair
$U(x,y)=(x,yx)$ and $V(x,y)=(xy^{-1},y)$, so $J(x,y)=(xy,yxy)$.
Evaluating $JU=UJ$ on $x$ gives $xy=xyx$, hence $x=1$, and on $y$
then gives $y^3=y^2$, hence $y=1$. The fiber dies with its boundary,
the meridian of $F$. Since $Y\setminus\nu T$ is simply connected by
Lemma~\ref{lem:exteriors}, this meridian normally generates
$\pi_1(Y\setminus\nu(F\cup T))$. A parallel
section puts the remaining bundle generators on the gluing boundary,
so $\pi_1(Z_g\setminus\nu T)=1$.

Since $V_g$ bounds over $B_0^\circ\x[0,1]$ and $\eu(V_g)=4(g-1)$,
we have $\sigma(Z_g)=0$ and $b_2^\pm(Z_g)=3+2g$. For spin $Y$ the
twists preserve the fiber spin structure whose
quadratic refinement satisfies $q(x_i)=q(y_i)=1$ and lift to its spin bundle
since $\pi_1(B_0^\circ)$ is free. After pullback a class from
$H^1(F;\Z_2)$ matches the boundary spin structures. For nonspin $Y$,
projection to $F$ and collapse of the normal circles along a collar
give a degree one map $Z_g\to Y$, preserving an odd square.
Freedman's theorem~\cite{Freedman} now gives the stated types for
$g=(m-3)/2$ and $g=(n-3)/2$, respectively. The sum is minimal since
$Y$ is minimal and $V_g$ aspherical~\cite{Usher}, and knot surgery
on $T$ gives both families as in Theorem~\ref{thm:main}.
\end{proof}

The two framings for surgery on a circle winding twice around the $S^1$ in $S^1\x S^3$
give a spin manifold $L_2$ and a nonspin manifold $L'_2$~\cite{MeierSpun}.
For $\pi_1=\Z_2$ and $b_2^+=b_2^-=m$, the three smoothable homeomorphism
types~\cite[Theorem~C]{HamKreck} are
\textup{(I)} $L_2\#(\#_m(\CP\#\CPb))$ with odd form,
\textup{(II)} $L_2\#(\#_m(S^2\x S^2))$ (spin), and
\textup{(III)} $L'_2\#(\#_m(S^2\x S^2))$ (nonspin with spin universal cover).

\begin{corollary}\label{cor:cyclic}
For every $m\ge3$, all three manifolds above admit infinitely many
pairwise nondiffeomorphic irreducible exotic copies. For odd $m$ these
copies are symplectic, and for even $m$ their universal covers are
symplectic.
\end{corollary}

\begin{proof}[Sketch of proof]
For odd $m$ take $(Y,T)$ as above with $b_2^\pm(Y)=m$ and simply
connected torus exterior. For \textup{(I)} and \textup{(II)}, sum $T$ with
the zero section of the torus bundle over $T^2$ with monodromies a quarter
turn and the identity~\cite[Lemma~4.1]{APSYtorus}. The sum preserves
$\eu,\sigma$, and killing the section meridian and base loops leaves the
fiber relation $x=y=x^{-1}$, giving $\pi_1=\Z_2$. Collapsing punctured
fibers to disks gives a degree one map to $Y$, preserving odd squares.
For spin $Y$ the bundle is $S^1$ times an oriented $3$-manifold and base
classes match the boundary spin structures.

For \textup{(III)} sum spin $Y$ with the symplectic bisection $S^1\x K$
in $S^1\x(S^1\x S^2)$, where $K$ is the $(2,1)$ curve on $S^1$ times
an equator. Killing its meridian and longitudes gives
$\Z^2/\gen{(1,0),(0,2)}=\Z_2$. Since $K$ winds twice, spin structures
extending over the bisection agree on its framed longitude. Turning the
normal disk once along $K$ adds the meridian, so one gluing fails to
match the spin structure from the $Y$ side. In the double base cover a
base class matches both sections along their primitive longitudes. All three sums are minimal symplectic and
irreducible~\cite{Usher,HK}.

For even $m$ use the double construction of~\cite[Lemmas~2.3--2.4]{BSS},
starting with $b_2^\pm(Y)=m-1$. In $\Sigma_2\x T^2$ write $u_i,v_i$ for
the handle pairs and $x,y$ for the torus generators. Sum $u_1\x x$ with
$T$ and do unit Luttinger surgeries on $v_1\x x'$, $u_2\x y$ along
$v_1,y$, with disjoint duals $u_1'\x y'$, $v_2'\x x''$. Off
$\Sigma=\Sigma_2\x\{\mathrm{pt}\}$ these operations kill $u_1,x,v_1,y$.
Double along $\Sigma\x S^1$ by $(p,z)\mapsto(jp,\pm z)$, where the deck involution $j$ of
$\Sigma\to N=T^2\#\mathbb{RP}^2$ exchanges handles, killing the remaining
generators. The isotopic gluings give the same simply connected minimal
symplectic double and irreducible quotients of $\eu=2m+2$, $\sigma=0$.
For nonspin $Y$, either an odd class orthogonal to $T$ survives or an
odd dual caps with a punctured $v_1'\x y'$ disjoint from $\Sigma$ and
the surgeries, giving \textup{(I)}.

For spin $Y$ the $+$ collar is $J\x S^1$, with $J$ the orientable
interval bundle over $N$. Spin structures on $J,T^2$ nonbounding on
$v_1,y$, respectively, give a product structure which survives surgery,
matches $Y$ and extends across this collar with bounding meridian,
giving \textup{(II)} (cf.~\cite[Theorem~3.1]{AMeven}). The $-$ collar
is the circle bundle of the complexified orientation line of $N$ pulled
back to $J$. Its mod two Euler class $w_1(N)^2$ is nonzero since $\eu(N)$
is odd, so by the Gysin sequence all spin structures agree on the fiber,
which is nonbounding in the
vertical framing. On the other side the meridian bounds, with spin
value fixed by the dual $\{\mathrm{pt}\}\x T^2$. This gives
\textup{(III)} with the same spin double.

For the families do fibered knot surgery on a parallel of $T$ before
taking the sums. The remaining torus $S^1\x U$, with $U$ a meridian of
the knot, still has simply connected exterior. In the link group, killing
the knot meridian leaves the unknot group of $U$, and the knot longitude
kills its generator by linking number one. The bounding Seifert longitude
preserves the spin gluings. The two lifted surgery tori are symplectic
up to orientation in separate copies of the $Y$ exterior, hence essential,
and the second remains symplectic after surgery on the first. Applying
the knot surgery formula~\cite[Theorem~3.1]{Chen} twice distinguishes
the covers and hence the quotients as in Theorem~\ref{thm:main}.
\end{proof}

\begin{remark}\label{rem:cyclic-extensions}
Larger sums of types \textup{(I)} and \textup{(II)} were already
realized in~\cite{BSS,AModd} and~\cite{AMeven}, respectively.
As in~\cite{BSS,LT}, one can also obtain irreducible exotic
families with $\pi_1=\Z_n$ for some $n>2$.
\end{remark}

\begin{remark}\label{rem:free}
Changing the torus surgeries gives free anti-symplectic involutions on exotic
copies of both manifolds in Theorem~\ref{thm:main}. The quotients belong
to infinite families of pairwise nondiffeomorphic irreducible exotic copies
of $L'_2\#(S^2\x S^2)$ and $L_2\#(\CP\#\CPb)$, respectively.
While these smallest examples of signature zero exotica with $b_2^+=1$ are fairly
interesting, we leave the details of their construction for later. At present we
have only a computer proof that their double covers are simply connected.
\end{remark}
\endgroup

\newpage
\section{Complex constructions}\label{sec:complex}

To prove Theorem~\ref{thm:complex}, we take the quotient of a complex fake projective plane $F$
by $\Z_3$, pass to an unramified cover of this complex orbifold, and
resolve the six $A_2$ singularities upstairs. We choose $F$ from the
Cartwright--Steger classes $(C18,3,\varnothing)$ for the spin example
and $(C10,2,\{17-\})$ for the nonspin example~\cite{CSplanes}.
JongHae Keum showed that an order three automorphism of an FPP has
three fixed points, with local action $(z_1,z_2)\mapsto
(\zeta z_1,\zeta^{-1}z_2)$ for a primitive cube root $\zeta$~\cite{Keum}.

Remove open cone neighborhoods of their images in $F/\Z_3$ to obtain
$W$ with three boundary lens spaces $L(3,2)$. The resolution replaces each cone by the
plumbing of two $(-2)$-spheres meeting once. In Appendix~\ref{app:covers}
we construct the following covers of $W$ by giving their monodromy.

\begin{lemma}\label{lem:cover-actions}
Each $W$ admits a connected $32$-fold cover $p\colon\widetilde W\to W$.
Over each boundary $L(3,2)$, two copies map \textup{1:1} and ten
$S^3$s map \textup{3:1}. The images of the six lens space groups normally generate
$\pi_1(\widetilde W)$, and $p^*\colon H^2(W;\Z_2)\to H^2(\widetilde W;\Z_2)$
is zero for $C18$ and injective for $C10$.
\end{lemma}

\begin{proof}[Proof of Theorem~\ref{thm:complex}]
Coning off $\partial\widetilde W$ extends $p$ to an orbifold cover
$Z\to F/\Z_3$. The cone on $S^3$ is a $4$-ball, and replacing
the six lens space cones by plumbings gives the resolution $Y$.
Since the plumbings are simply
connected and the six lens space groups normally generate
$\pi_1(\widetilde W)$, van Kampen gives $\pi_1(Y)=1$.

Let $d=\deg p$ and let $s$ count the lens spaces upstairs. A common
smooth finite cover with the ball quotient $F$ shows that $Z$ is projective
with ample orbifold canonical class of square $dK_F^2/3$. As in Keum's proof
for $F/\Z_3$~\cite{Keum}, this class pulls back to $K_Y$,
so $Y$ is minimal of general type.

Since $F$ has the Betti numbers of $\CP$, the \textup{3:1} cover by $F$ minus
three balls gives $\eu(W)=0$. Of the $3d$ boundary sheets, $s$ give lens
spaces and the rest give $(3d-s)/3$ copies of $S^3$. Each plumbing adds three
to $\eu(Y)$ and each ball one. Using $K_F^2=9$ and $3\sigma(Y)=K_Y^2-2\eu(Y)$, we get
\[
 K_Y^2=3d,\qquad \eu(Y)=3s+(3d-s)/3,\qquad 9\sigma(Y)=3d-16s .
\]
Here $d=32$, $s=6$, giving $K_Y^2=96$, $\eu(Y)=48$, $\sigma(Y)=0$
and $b_2^\pm(Y)=23$.

To determine parity it suffices to check $\widetilde W$, since the balls
and plumbings are spin and each boundary has a unique spin structure.
As $K_F^2$ is odd, $F$ is nonspin, and so is $F$ minus three balls.
A spin structure on $W$ would lift to this \textup{3:1} cover, so $w_2(W)\ne0$.
By Lemma~\ref{lem:cover-actions}, $w_2(\widetilde W)=p^*w_2(W)$
vanishes for $C18$ and survives for $C10$, giving the spin and nonspin
examples, respectively. Freedman's theorem identifies their
homeomorphism types, and Taubes's nonvanishing distinguishes them
from the standard sums as in Section~\ref{sec:families}.
\end{proof}

Here $\sigma=0$ forces $d=16k$, $s=3k$, giving $b_2^\pm=12k-1$ when
$b_1=0$. The $C18$ deck involution on $\widetilde W$ extends over the cones
and lifts holomorphically to the minimal resolution. In fact, the lift
is free since the involution exchanges boundary components in pairs
and hence also the corresponding balls and plumbings. Its quotient is a compact complex surface
with $\pi_1=\Z_2$ and $b_2^\pm=11$. The final paragraph of Appendix~\ref{app:covers}
describes this involution on the sheets and identifies the quotient as
the nonspin type with spin universal cover in Corollary~\ref{cor:cyclic}.

\begin{remark}[Cannot reach three summands this way]\label{rem:fpp-obstruction}
For any FPP $F$ and any group $G$ of holomorphic automorphisms of $F$,
no finite orbifold cover $Z\to F/G$ has minimal resolution $Y$ with
$b_1(Y)=0$ and $b_2^+(Y)=b_2^-(Y)=3$. Thus neither homeomorphism type
in Theorem~\ref{thm:main} arises this way.

Keum described the singularities and minimal resolutions of the
quotients $F/G$~\cite{Keum}. He also showed, for example, that a rational
homology projective plane with just three $A_2$ singularities is a
quotient of an FPP by $\Z_3$ if its minimal resolution is minimal of
general type with $p_g=0$ and $K^2=3$~\cite{Keum12}. For our count,
the nontrivial point stabilizers have prime order three or seven, so each singularity
upstairs has the type of its image. Let $a,b$ count the two types on $Z$. An order three point has weights
$(1,2)$ and a $(-2,-2)$ resolution chain, replacing the orbifold Euler
contribution $1/3$ by $3$ and leaving $K^2$ unchanged. For order
seven the weights are $(1,3)$ and the chain is $(-3,-2,-2)$, replacing
$1/7$ by $4$ and lowering $K^2$ by $3/7$ as in
Keum's calculation~\cite{Keum}.
Since $K_Z^2=3\eu_{\mathrm{orb}}(Z)$, the proposed Betti numbers give
$8a+12b=3\eu(Y)-K_Y^2=\eu(Y)-3\sigma(Y)=8$.

Thus $a=1$, $b=0$. Over the image of this single point one sheet is
fixed and the rest come in triples, so the degree is $1\bmod 3$. Every
order three point downstairs must then have a fixed sheet, and Keum's
list gives at least three such points whenever there is one.
\end{remark}

\begin{remark}[A small complex surface of positive signature]\label{rem:positive}
A similar construction provides the following interesting complex surface
of signature $+1$. Take the FPP with group \texttt{index3aFP} in the
class $(15,2,\varnothing)$~\cite{CSplanes} and let the generators $A,B$ of $\pi_1(W)$ act on $\Z_{19}$
by $i\mapsto i+1$ and $i\mapsto8i$. The three boundary generators
fix $17,4,15$, respectively, with all other sheets in triples. Thus
$(d,s)=(19,3)$ gives a minimal surface of general type with
$(K_Y^2,\eu(Y),\sigma(Y))=(57,27,1)$. Killing the three lens space
groups gives $\pi_1(Y)=\gen{u,v\mid u^2,v^2,(uv)^2}\cong(\Z_2)^2$.
Since this group is finite, $b_1(Y)=0$. With $b_2^+(Y)=(27-2+1)/2=13$,
this is by far the smallest compact complex surface of general type with
positive signature and finite fundamental group that the author is aware of.
\end{remark}

\linespread{1}\selectfont
\appendix
\section{The two special complex covers}\label{app:covers}

To construct the two covers in Lemma~\ref{lem:cover-actions}, we label
the sheets from $1$ to $32$ and record where a lift of each loop ends.
Permutations satisfying the relators of $\pi_1(W)$ specify a cover,
which is connected if the action is transitive. For a generator of a
boundary lens space group, a fixed sheet gives a copy of $L(3,2)$ to
which we attach a plumbing, and a cycle of length three gives an $S^3$ to which we
attach a ball. We then check that the six lens space groups normally
generate $\pi_1(\widetilde W)$, so the resolution is simply connected,
and determine whether $w_2(W)$ survives upstairs, which distinguishes
the spin and nonspin examples.

\begin{proof}[Proof of Lemma~\ref{lem:cover-actions}]
Write $\Gamma=\pi_1(W)$ and let $\Pi=\pi_1(F)$ be its subgroup of
index three. Denote the images of generators of the three boundary
lens space groups by $t_1,t_2,t_3$. In the catalogue
words $\bar g=g^{-1}$, and permutations act on the right, from left to right.
For $C18$, Cartwright and Steger's \texttt{index3bFP} and
\texttt{index9aFP} denote $\Gamma$ and $\Pi$ of indices $3$ and $9$
in their group on $A,B,Z$. Setting $h=(BZ)^{-1}$, take
$x=B^{-1}Z^{-1}$, $y=ZAZ^{-1}$ and $z=hAhZAhZ^{-1}$ for $\Gamma$.
Their permutations are
\[
\begin{aligned}
\rho(x)&=(1,2,3)(4,5,6)\cdots(28,29,30), \\
\rho(y)&=(1,21,19,3)(2,31,5,30,24,7,13,28)(4,22,17,14) \\
 &\phantom{=}(6,11,26,23)(8,12,29,18,9,20,32,25)(10,15,16,27), \\
\rho(z)&=(1,2,11,5,17,24,4,13)(3,31,21,32)(6,30,22,7)(8,10,29,26) \\
 &\phantom{=}(9,14,28,16)(12,27,18,19,20,23,25,15).
\end{aligned}
\]
In the cyclic quotient by $\Gamma$ we have $A=1$, $B=Z^{-1}$, so only
$BZ$ of the four catalogue torsion representatives lifts closed. The
coset representatives $1,Z,Z^{-1}$ give the conjugates $g(BZ)g^{-1}$,
one per boundary component, whose cubes vanish by $(BZ)^3=1$.
Rewriting these in order in $x,y,z$ gives $t_1=v$, $t_2=\bar{x}$ and
$t_3=\bar{u}\bar{z}\bar{x}vyvzuxyx\bar{z}xzuxyxv$, where
$u=\bar{y}x\bar{z}y$ and $v=yxzu$. The relators are
\[
\begin{gathered}
x^3,\quad \bar{y}zx\bar{z}uyx\bar{z}y^3xz\bar{x},\quad yxzxz\bar{x}\bar{u}x\bar{z}v,\quad yxz\bar{y}^2\bar{u}z\bar{x}\bar{y}\bar{x}\bar{u}\bar{z}^2, \\
z\bar{x}ux\bar{z}\bar{x}uv,\quad \bar{x}\bar{y}\bar{x}\bar{v}\bar{y}\bar{u}zxyx\bar{z},\quad zxz\bar{x}v\bar{x}\bar{z}yxyxz\bar{x}\bar{y}x\bar{y}, \\
\bar{x}\bar{z}\bar{x}\bar{y}^2\bar{u}z\bar{x}\bar{y}x\bar{y}\bar{x}\bar{u}\bar{z}\bar{x}z,\quad \bar{x}\bar{y}\bar{x}\bar{u}\bar{z}x\bar{y}zx\bar{u}zxyx\bar{z}^2,\quad yv\bar{z}xzuxy\bar{x}zxy^2x\bar{z}, \\
(\bar{y}\bar{x}\bar{u}\bar{z}x)^3,\quad yx\bar{y}^2\bar{u}z^3\bar{x}\bar{y}\bar{x}\bar{z}^3\bar{v}x,\quad z^2x\bar{z}\bar{y}^2\bar{u}z\bar{x}(\bar{z}\bar{x}\bar{z})^2yx\bar{z}xzuxyx\bar{y}\bar{v}.
\end{gathered}
\]

For $C10$, $\Gamma$ is the catalogue group on $(x,y,z)=(A,B,C)$ and
$\Pi$ is \texttt{index3FP}. Here
\[
\begin{aligned}
\rho(x)&=(1,2,\dots,14)(15,16,\dots,28)(29,30)(31,32), \\
\rho(y)&=(1,21,12,25,3,30,20,9,5,6,15,29) \\
 &\phantom{=}(2,24,18,22,32,27,28,17,13,10,31,8)(4,14) \\
 &\phantom{=}(7,23,26,11,19,16), \\
\rho(z)&=(1,28,8,21)(2,20,9,27)(3,18,22,13)(4,31,16,29)(5,7) \\
 &\phantom{=}(6,10,25,15)(11,32,23,30)(12,14)(17,19)(24,26).
\end{aligned}
\]
Here the catalogue's three finite order representatives
(\texttt{FOList}) give the boundary words
$t_1=\bar{z}y$, $t_2=\bar{y}x$ and $t_3=\bar{y}\bar{x}y^2\bar{z}y\bar{x}$.
The relators are
\[
\begin{gathered}
(y\bar{x})^3,\quad (y\bar{z})^3,\quad \bar{z}y\bar{x}^2yz\bar{y}x^2\bar{y},\quad yx\bar{y}z\bar{y}\bar{x}^2y^2\bar{z}y^2\bar{z},\quad (\bar{x}y\bar{x})^2\bar{x}zx^3\bar{y}z\bar{y}, \\
xy^2x^3\bar{y}z^2\bar{y}\bar{z}y\bar{z}\bar{y}z,\quad \bar{y}x^2\bar{y}(\bar{y}x^2\bar{y}z)^2z^3,\quad yx\bar{z}y\bar{z}\bar{y}zxy\bar{x}yz\bar{y}z\bar{x}^3y\bar{z}, \\
xyx\bar{y}z\bar{y}(\bar{y}xyx\bar{z}y\bar{z})^2\bar{y},\quad x\bar{y}zx^2\bar{y}^2\bar{x}y\bar{x}^2zxy^3\bar{x}^2z\bar{y}x.
\end{gathered}
\]

The permutations satisfy the relators and act transitively, giving a
connected cover with group $H=\operatorname{Stab}_\Gamma(1)$. Each $t_j$
fixes the two sheets listed below and permutes the other thirty in
triples, giving the claimed boundary components.
To check normal generation, attach disks
along the six fixed lifts in
the lifted presentation complex. Collapse the spanning tree obtained
from sheet~$1$ in discovery order along $x,\bar x,y,\bar y,z,\bar z$.
The resulting presentations on $96-31=65$ generators and $422$ or $326$
relators reduce to the empty one by solving for letters occurring once
in a relator and shortening words with cyclic conjugates of other relators.

For parity we work over $\Z_2$. The exponent sum matrices have rank
three, so $H^1(W)=0$. The three boundary homology spheres give $b_3(W)=2$
by duality and the pair sequence, hence $b_2(W)=1$ since $\eu(W)=0$.
The universal cover of $W$ is the open unit ball in $\mathbb C^2$
with a locally finite family of open $4$-balls removed. It is therefore
$2$-connected, so group cohomology computes the required pullback in degree two.

Label the base relators by a vector $s$. The matrix $P$ repeats each
label on its lifts, and $D$ counts edges in the lifted relators modulo
two. Thus $s$ pulls back to a coboundary precisely when
$Ps\in\operatorname{im}D$. Every such $s$ is a cocycle on $K(\Gamma,1)$,
since the attaching spheres of the $3$-cells lift. Subtracting the
three independent coboundaries downstairs gives kernel dimension
$n+\operatorname{rk}D-\operatorname{rk}[D\ P]-3$, where $n$ is the
number of base relators.
\[
\begin{array}{c|ccc|cc}
 &\operatorname{Fix}(t_1)&\operatorname{Fix}(t_2)&\operatorname{Fix}(t_3)
 &\operatorname{rk}D&\operatorname{rk}[D\ P]\\ \hline
 C18 &\{11,23\}&\{31,32\}&\{5,25\}&65&74\\
 C10 &\{9,22\}&\{6,28\}&\{27,32\}&65&72
\end{array}
\]
The last two columns give $13+65-74-3=1$ for $C18$ and $10+65-72-3=0$ for $C10$,
as claimed.
\end{proof}

It remains to verify the additional claim in Section~\ref{sec:complex}
that the spin $C18$ surface has a free involution with nonspin quotient.
The orbit of $\{31,32\}$ under the $C18$ permutations partitions the
$32$ sheets into sixteen pairs. Swapping partners commutes with monodromy,
so it gives a deck involution. For each boundary generator it exchanges
the two fixed sheets and pairs the three-cycles, hence it exchanges the
attached balls and plumbings and acts freely on the resolution.
On the sixteen pairs the ranks of $D$ and $[D\ P]$ are $32$ and $42$,
so the pullback on $H^2(W;\Z_2)$ has kernel of dimension
$13+32-42-3=0$. Thus $w_2(W)$ survives on the degree $16$ cover and
the quotient surface is nonspin.

\smallskip
\bibliographystyle{amsplain}
\bibliography{3S2xS2}

\end{document}